# RUIN PROBABILITIES AND DECOMPOSITIONS FOR GENERAL PERTURBED RISK PROCESSES[1]

By Miljenko Huzak,[2] Mihael Perman, Hrvoje Šikić[3]
and Zoran Vondraček[4]

*University of Zagreb, University of Ljubljana and University of Zagreb*

We study a general perturbed risk process with cumulative claims modelled by a subordinator with finite expectation, with the perturbation being a spectrally negative Lévy process with zero expectation. We derive a Pollaczek–Hinchin type formula for the survival probability of that risk process, and give an interpretation of the formula based on the decomposition of the dual risk process at modified ladder epochs.

**1. Introduction.** The classical Cramér–Lundberg model in insurance assumes that the risk process $(R(t), t \geq 0)$ is given by $R(t) = ct - \sum_{i=1}^{N(t)} Y_i$, where $c > 0$ is the premium rate, $(Y_i, i \in \mathbb{N})$ is an i.i.d. sequence of nonnegative random variables modelling individual claims, and $(N(t), t \geq 0)$ is a homogeneous Poisson process of rate $\lambda > 0$, independent of $(Y_i, i \in \mathbb{N})$. Hence the cumulative claim process is modelled by the compound Poisson process $\sum_{i=1}^{N(t)} Y_i$. Let $F$ denote the distribution function of $Y_i$, and let $\mu = \mathbb{E}Y_i$. The central question for the model is the computation of the ruin probability in infinite time, given initial capital $x > 0$, defined by

$$\vartheta(x) := \mathbb{P}(R(t) + x < 0 \text{ for some } t > 0).$$

In case $c \leq \lambda\mu$, this quantity is identically equal to 1. Hence, one always assumes the net profit condition $c > \lambda\mu$, and defines the parameter $\rho :=$

Received March 2003; revised June 2003.
[1]Supported by a bilateral Croatian–Slovenian research grant.
[2]Supported in part by MZT Grant 0037115 of the Republic of Croatia.
[3]Supported in part by MZT Grant 0037118 of the Republic of Croatia.
[4]Supported in part by MZT Grant 0037107 of the Republic of Croatia.
*AMS 2000 subject classifications.* Primary 60J25; secondary 60G51, 60J75, 60K30, 91B30.
*Key words and phrases.* Risk theory, ruin probability, Pollaczek–Hinchin formula, subordinator, spectrally negative Lévy process.

$\lambda\mu/c < 1$. Instead of studying the ruin probability, one can equivalently consider the survival probability $\theta(x) := 1 - \vartheta(x)$, which is more convenient. One of the few explicit results for the survival probability is the Pollaczek–Hinchin formula:

$$(1.1) \qquad \theta(x) = (1-\rho) \sum_{n=0}^{\infty} \rho^n F_I^{n*}(x),$$

where $F_I(x) = (1/\mu) \int_0^x (1 - F(t))\, dt$ is the integrated tail distribution. Formula (1.1) is usually derived via renewal arguments. The resulting integro-differential equation for $\vartheta$ is solved using Laplace transforms. The explanation of the formula is given by considering the supremum of the dual process $\widehat{R}(t) := -R(t)$. By the net profit condition $\widehat{R}(t)$ drifts to $-\infty$, hence the supremum is a.s. finite, and clearly, $\theta(x) = \mathbb{P}(\sup_{0 \leq t < \infty} \widehat{R}(t) \leq x)$. It is easy to see that $\sup_{0 \leq t < \infty} \widehat{R}(t)$ is a sum of geometrically many i.i.d. random variables. It is not, however, quite as easy to determine the distribution of these variables. Usually fluctuation theory is used. We refer the reader to [1] and [9] for details.

In this paper we are interested in generalizations of the Cramér–Lundberg model, which lead to the same type of the Pollaczek–Hinchin formula for the survival probability, and which admit an explanation of the formula by decomposition of the supremum of the dual process in the random sum of ladder heights. One possible generalization of the model is to allow for additional uncertainties in the cumulative claims and/or in the premium income. These uncertainties may be the result of fluctuations in the claim arrival intensity, the number of insurees, inflation or surplus investment (see [9], page 568). Dufresne and Gerber [3] considered the risk process $(R(t), t \geq 0)$ perturbed by a multiple of standard Brownian motion $(W(t), t \geq 0)$, and defined $X(t) := R(t) + \varsigma W(t)$, $\varsigma > 0$. Using renewal arguments, they derived the formula

$$(1.2) \qquad \theta(x) = (1-\rho) \sum_{n=0}^{\infty} \rho^n (G^{(n+1)*} * F_I^{n*})(x).$$

The parameter $\rho$ and the distribution function $F_I$ are the same as in the unperturbed model, while $G$ is an exponential distribution function with parameter $2c/\varsigma^2$. They also gave the following interpretation of the formula (1.2): Let $\sigma_1, \sigma_2, \ldots$ be the moments when a new supremum of the dual process $\widehat{X}(t) := -X(t)$ is reached by a jump of the claim process $\sum_{i=1}^{N(t)} Y_i$. Then the number of such moments has geometric distribution with parameter $\rho$, $G$ is the distribution function of the supremum of $\widehat{X}(t)$ just before $\sigma_1$, and $F_I$ is the conditional distribution of the overshoot over the previous supremum, given $\sigma_1 < \infty$.

Furrer [5] considered the process $X(t) = R(t) + Z_\alpha(t)$, where $R$ is the classical risk process, and $Z_\alpha$ is an $\alpha$-stable Lévy process with no positive jumps, $1 < \alpha < 2$. He used the explicit formula for the Laplace exponent of the infimum of $X(t)$ due to Zolotarev [13] to obtain formula (1.2) for the survival probability of $X(t)$. The distribution function $G$ is explicitly identified as the Mittag–Leffler distribution given by $1 - G(x) = \sum_{n=0}^{\infty} (-cx^{\alpha-1})^n / \Gamma(1 + (\alpha - 1)n)$.

Schmidli [11] gives a nice interpretation of $G$ as the distribution of the supremum of the dual process $\widehat{X}$ just before the first time the process $\widehat{X}$ exceeds its previous supremum by a jump of the cumulative claim process. His setting is more general in the sense that the cumulative claim process is generated by a stationary, ergodic, marked point process.

Another possible generalization of the classical risk process is to allow a different cumulative claim process. Dufresne, Gerber and Shiu [4] considered the model in which the claim process was modelled by a Gamma process. Such a process has infinitely many jumps in finite intervals. Yang and Zhang [12] studied this model perturbed by a Brownian motion. Using the approach in [5], they derived a formula of the type (1.2) with $F_I$ replaced by an exponential integral type distribution, and $G$ is again the exponential distribution.

The goal of this paper is twofold: (1) to extend the Pollaczek–Hinchin formula to the more general setting of the spectrally negative Lévy processes while retaining the risk insurance theory interpretation; (2) to give a unified and transparent approach to the problem by use of well-developed and powerful fluctuation theory for Lévy processes. We will consider a general perturbed risk process $X(t) = ct - C(t) + Z(t)$, where $(C(t), t \geq 0)$ is a cumulative claim process, and $(Z(t), t \geq 0)$ is a perturbation. Note that the cumulative claim process has to be increasing. Therefore, if one wants to stay in the realm of processes with stationary independent increments, the only choice for modelling $(C(t), t \geq 0)$ is subordinators. Hence, we assume that $(C(t), t \geq 0)$ is a subordinator (without drift) having finite expectation satisfying the net profit condition $c - \mathbb{E}C(1) > 0$. The perturbation is modelled by a Lévy process $(Z(t), t \geq 0)$ with no positive jumps, having zero expectation. The assumption that the expectation is zero is inconsequential, since $\mathbb{E}Z(1)$ can always be moved to the premium rate. In the analysis of the risk process $(X(t), t \geq 0)$, we will rely heavily on fluctuation theory for general Lévy processes, which is particularly explicit for processes with no positive jumps. For background on these results, we refer the reader to the book by Bertoin [2].

Our first result is the formula for the survival probability for the process $X$ which is proved in Section 3:

$$\begin{aligned}(1.3)\qquad \theta(x) &:= \mathbb{P}\bigg(\inf_{0\leq t<\infty} X(t) > -x\bigg)\\ &= (1-\rho)\sum_{n=0}^{\infty} \rho^n (G^{(n+1)*} * H^{n*})(x).\end{aligned}$$

We essentially follow the approach from [5], and obtain explicitly the parameter $\rho$ and the distribution functions $G$ and $H$ appearing in the formula. It turns out that $G$ can be identified as the distribution function of the absolute supremum of the process $(-ct - Z(t), t \geq 0)$, while $H$ is related to the subordinator $C(t)$ only, and can be thought of as the integrated tail distribution of jumps. In Section 4 we give an interpretation of formula (1.3) by decomposing the supremum of the dual process $\widehat{X}(t) := -X(t)$ into the random sum of modified ladder heights. In order to do this, we first show that the times when the new supremum of $\widehat{X}(t)$ is reached by a jump of the subordinator are discrete. Let $0 < \sigma_1 < \sigma_2 < \cdots$ be those times, and let $\tilde{G}$ be the distribution function of $\widehat{S}(\sigma_1-)$, where $\widehat{S}(t) := \sup_{0\leq s\leq t} \widehat{X}(s)$. We show that $\widehat{S}(\sigma_1-)$ and the overshoot $\widehat{S}(\sigma_1) - \widehat{S}(\sigma_1-)$ are conditionally independent given $\sigma_1 < \infty$, and identify the conditional distribution of the overshoot with $H$. Using the strong Markov property at times $\sigma_i$, we rederive formula (1.3) with $\tilde{G}$ instead of $G$ (and the same $\rho$). This clearly implies that $\tilde{G} = G$, yielding the required interpretation. Our results are more general and cover the results obtained in [5], [12] and [11] (in Lévy case).

Another interpretation of formula (1.3) is provided by looking at the ladder height process of $\widehat{X}$. The ladder height process is obtained by time-changing $\widehat{S}(t)$ by the inverse local time at zero of the reflected process $\widehat{S}(t) - \widehat{X}(t)$. This process records only values where the new supremum is reached, and consequently, contains all the relevant information on the distribution of the supremum of $\widehat{X}(t)$. In Section 5 the results of Section 4 are reinterpreted and improved in terms of the ladder height process.

We end this Introduction by noting that in a very recent paper Klüppelberg, Kyprianou and Maller [7] study ruin probabilities for general Lévy insurance risk process (not necessarily spectrally negative) drifting to $-\infty$. They are mostly concerned with the asymptotic results for the first passage time and overshoot behavior at high levels.

**2. Setting and notation.** Let $(\Omega, \mathcal{F}, \mathbb{P})$ be a probability space on which all random variables will be defined. As explained in the Introduction, we model the cumulative claim process by a subordinator $C = (C(t), t \geq 0)$ without a drift. Let $\nu$ be the Lévy measure of $C$; that is, $\nu$ is a $\sigma$-finite measure on $(0,\infty)$ satisfying $\int_{(0,\infty)}(x\wedge 1)\nu(dx)<\infty$. The Laplace exponent of $C$ is defined by

$$\Phi_C(\beta):=\int_{(0,\infty)}(1-e^{-\beta x})\nu(dx)$$

so that

$$\mathbb{E}[\exp\{-\beta C(t)\}]=\exp\{-t\Phi_C(\beta)\}.$$

Note that

$$\mathbb{E}C(1)=\Phi'_C(0+)=\int_{(0,\infty)}x\nu(dx)=\int_0^\infty \nu(x,\infty)\,dx,$$

where the last equality follows by integration by parts. As explained in the Introduction, we assume throughout that $\mathbb{E}C(1)<\infty$. Let

(2.1) $$H(x):=\frac{1}{\mathbb{E}C(1)}\int_0^x \nu(y,\infty)\,dy.$$

Then $H$ is an absolutely continuous distribution function with density $h(x)=\nu(x,\infty)/\mathbb{E}C(1)$. We call $H$ the integrated tail distribution. The Laplace transform of $H$ is given by

(2.2) $$\mathcal{L}H(\beta):=\int_0^\infty e^{-\beta x}H(dx)=\int_0^\infty e^{-\beta x}h(x)\,dx=\frac{1}{\mathbb{E}C(1)}\frac{\Phi_C(\beta)}{\beta}.$$

Let $\Delta C(t)=C(t)-C(t-)$. It is well known that $(\Delta C(t)\colon t\geq 0)$ is a Poisson point process with characteristic measure $\nu$ and state space $(0,\infty)\cup\{\partial\}$. The cemetery state $\partial$ is added to keep with the standard definition of a Poisson point process (cf. [8], page 435). It is assumed that the process is in state $\partial$ whenever there is no jump and $\nu(\{\partial\})=0$. Moreover, one has $C(t)=\sum_{0<s\leq t}\Delta C(s)$.

We model the risk process $R=(R(t),\,t\geq 0)$ as $R(t)=ct-C(t)$, where $c>0$ is the premium rate. Then $R$ is a Lévy process with no positive jumps (i.e., spectrally negative Lévy process). The Laplace exponent $\psi_R$ of $R$ is defined by relation

$$\mathbb{E}[\exp\{\beta R(t)\}]=\exp\{t\psi_R(\beta)\}.$$

Clearly, $\psi_R(\beta)=c\beta-\Phi_C(\beta)$. It is important to note that $R(t)$ stays positive in a neighborhood of $t=0$, implying that ruin (with zero initial capital) does not occur immediately. This follows from the fact that $\lim_{t\to 0}C(t)/t=0$ ([2], page 84 or 192). From now on we assume that the net profit condition $c>\mathbb{E}C(1)$ holds, and let $d:=c-\mathbb{E}C(1)$. It follows that $\mathbb{E}R(1)=\psi'_R(0+)=d>0$, which implies that $R$ drifts to $+\infty$. We also introduce the parameter $\rho:=\mathbb{E}C(1)/c\in(0,1)$.

The perturbation $Z = (Z(t), t \geq 0)$ of the risk process $R$ will be modelled by a spectrally negative, mean zero, Lévy process. Its Lévy measure $\Pi_Z$ is an infinite $\sigma$-finite measure on $(-\infty, 0)$ satisfying the usual condition

$$\text{(2.3)} \quad \int_{(-\infty,0)} (x^2 \wedge 1) \Pi_Z(dx) < \infty,$$

and the additional condition

$$\text{(2.4)} \quad \int_{(-\infty,-1)} |x| \Pi_Z(dx) < \infty,$$

which ensures finite expectation of $Z$. The Laplace exponent of $Z$ is given by

$$\text{(2.5)} \quad \psi_Z(\beta) := \frac{\varsigma^2}{2} \beta^2 + \int_{(-\infty,0)} (e^{\beta x} - 1 - \beta x) \Pi_Z(dx),$$

where $\varsigma \geq 0$, and integrability of the integrand follows from condition (2.4). Further, $\mathbb{E}Z(1) = \psi_Z'(0+) = 0$ (e.g., [10], page 163). Note that we allow $Z$ to be identically zero (both $\Pi_Z = 0$ and $\varsigma = 0$). However, $Z$ cannot be compound Poisson because such processes cannot have $\mathbb{E}Z(t) = 0$. Let us point out that our setting includes the Brownian perturbation ($\varsigma > 0$, $\Pi_Z = 0$), and also the perturbation by $\alpha$-stable spectrally negative Lévy process for $\alpha \in (1, 2)$ ($\varsigma = 0$, $\Pi_Z(dx) = (a/|x|^{\alpha+1}) \mathbb{1}_{(-\infty,0)} dx$).

Finally, we define the general perturbed risk process $X = (X(t), t \geq 0)$ as

$$X(t) := R(t) + Z(t) = ct - C(t) + Z(t),$$

where $C$ and $Z$ are independent processes. The process $X$ is a spectrally negative Lévy process with finite positive expectation $\mathbb{E}X(1) = c - \mathbb{E}C(1) = d > 0$. Therefore, $\lim_{t \to \infty} X(t) = +\infty$ a.s., that is, $X$ drifts to infinity. Let $\mathcal{F}^0(t) := \sigma(C(s), Z(s), 0 \leq s \leq t)$, and let $\mathcal{F} = (\mathcal{F}(t), t \geq 0)$ be the filtration obtained in the usual way by augmenting $\mathcal{F}^0(t)$. Clearly, $X(t)$ is $\mathcal{F}(t)$-measurable. The Laplace exponent $\psi$ of $X$, defined by the relation

$$\mathbb{E}[\exp\{\beta X(t)\}] = \exp\{t \psi(\beta)\},$$

is, due to independence of $C$ and $Z$, given by

$$\psi(\beta) = c\beta - \Phi_C(\beta) + \psi_Z(\beta), \qquad \beta \geq 0.$$

Since $\psi$ is strictly convex and $\psi'(0+) = \mathbb{E}X(1) > 0$, $\psi$ is strictly increasing on $[0, \infty)$, and therefore has a strictly increasing inverse $\Phi \colon [0, \infty) \to [0, \infty)$. Since $\psi(0) = 0$, it follows that $\Phi(0) = 0$.

In the sequel, we will be interested in the function $\theta \colon [0, \infty) \to [0, 1]$ defined by

$$\text{(2.6)} \quad \theta(x) := \mathbb{P}(X(t) \geq -x, \text{ for all } t \geq 0).$$

This function is the survival probability of the general perturbed risk process $X$ starting with the initial capital $x \geq 0$. The initial behavior of $X$ determines $\theta$ at zero. If there is no perturbation, that is, if $X = R$, then, as said before, $X$ remains positive (a.s.) for an initial period of time, and hence $\theta(0) > 0$. On the other hand, if $Z \neq 0$, then $X$ is of unbounded variation, hence the point 0 is regular for $(-\infty, 0)$ ([2], page 192). Thus $X$ hits the interval $(-\infty, 0)$ immediately, implying $\theta(0) = 0$.

**3. Laplace transform approach.** In this section we derive the Pollaczek–Hinchin formula for the survival probability using the explicit form of the Laplace transform of the absolute infimum of $X$. Let $I(\infty) := \inf_{0 \leq s < \infty} X(s)$ and $I(t) := \inf_{0 \leq s \leq t} X(s)$. The fluctuation theory for Lévy processes provides the following formula for the Laplace transform of the infimum evaluated at an independent exponential time $\tau(q)$ with parameter $q > 0$ (see [2], page 192):

$$\mathbb{E}[\exp\{\beta I(\tau(q))\}] = \frac{q(\Phi(q) - \beta)}{\Phi(q)(q - \psi(\beta))}, \qquad \beta > 0.$$

Letting $q \downarrow 0$, and using $I(\tau(q)) \xrightarrow{\mathbb{P}} I(\infty)$, it follows that

$$(3.1) \qquad \mathbb{E}[\exp\{\beta I(\infty)\}] = \psi'(0+)\frac{\beta}{\psi(\beta)} = d\frac{\beta}{\psi(\beta)}, \qquad \beta > 0.$$

Let us introduce for a moment the following notation: $Y(t) = ct + Z(t)$ and $\psi_Y(\beta) = c\beta + \psi_Z(\beta)$. By the same argument as above it follows that

$$(3.2) \quad \mathbb{E}\left[\exp\left\{-\beta\left(-\inf_{0 \leq t < \infty} Y(t)\right)\right\}\right] = \psi'_Y(0+)\frac{\beta}{\psi_Y(\beta)} = c\frac{\beta}{\psi_Y(\beta)}, \qquad \beta > 0.$$

Let $G$ denote the distribution function of $-\inf_{0 \leq t < \infty} Y(t) = \sup_{0 \leq t < \infty}(-ct - Z(t))$. Then the last formula says that

$$(3.3) \qquad \mathcal{L}G(\beta) := \int_0^\infty e^{-\beta x} G(dx) = c\frac{\beta}{\psi_Y(\beta)}, \qquad \beta > 0.$$

Recall formulae (2.1) and (2.2) from Section 2:

$$H(x) := (1/\mathbb{E}C(1))\int_0^x \nu(y, \infty)\,dy \quad \text{and} \quad \mathcal{L}H(\beta) = \Phi_C(\beta)/(\mathbb{E}C(1)\beta).$$

Also recall that $\rho = \mathbb{E}C(1)/c$, hence $d/c = (c - \mathbb{E}C(1))/c = 1 - \rho$. Now we compute $d\beta/\psi(\beta)$ in terms of $\rho$, $\mathcal{L}G$ and $\mathcal{L}H$. This idea comes from [5]:

$$d\frac{\beta}{\psi(\beta)} = d\frac{1}{\psi_Y(\beta)/\beta - \Phi_C(\beta)/\beta}$$

$$= d\frac{1}{c/\mathcal{L}G(\beta) - \mathbb{E}C(1)\mathcal{L}H(\beta)}$$

$$= \frac{d}{c} \frac{\mathcal{L}G(\beta)}{1 - \rho \mathcal{L}G(\beta)\mathcal{L}H(\beta)}$$

$$= (1-\rho)\mathcal{L}G(\beta) \sum_{n=0}^{\infty} (\rho \mathcal{L}G(\beta)\mathcal{L}H(\beta))^n.$$

By inverting the Laplace transform, we obtain the following theorem.

THEOREM 3.1. *The survival probability of the general perturbed risk process $X$ is given by*

$$\theta(x) = \mathbb{P}(I(\infty) \geq -x)$$
(3.4)
$$= (1-\rho) \sum_{n=0}^{\infty} \rho^n (G^{(n+1)*} * H^{n*})(x), \qquad x \geq 0.$$

We point out that $H$ depends only on the subordinator $C$, while $G$ depends on the premium rate $c$ and the perturbation $Z$. Brownian perturbations were considered in [5] and $\alpha$-stable ones in [3] and [12]. In both cases the distribution $G$ is given explicitly. If there is no perturbation, $Z = 0$, then $\mathcal{L}G(\beta) = 1$, and consequently, the distribution function $G$ can be omitted from formula (3.4).

**4. Decomposition of the supremum of $\widehat{X}$.** Let $\widehat{X}(t) := -X(t) = -ct + C(t) - Z(t)$ denote the dual process of $X$. Let

$$\widehat{S}(t) := \sup_{0 \leq s \leq t} \widehat{X}(s) \quad \text{and} \quad \widehat{S}(\infty) := \sup_{0 \leq s < \infty} \widehat{X}(s).$$

Since $\widehat{X}$ drifts to $-\infty$, $\widehat{S}(\infty) < \infty$ a.s. Introduce the following notation: $\widehat{I}(t) := \inf_{0 \leq s \leq t} \widehat{X}(s)$ and $S(t) := \sup_{0 \leq s \leq t} X(s)$. Clearly, $-\widehat{I}(t) = S(t)$. By a time reversal argument, $-\widehat{I}(t) \stackrel{d}{=} \widehat{S}(t) - \widehat{X}(t)$, and hence

(4.1) $$\widehat{S}(t) - \widehat{X}(t) \stackrel{d}{=} S(t).$$

In this section we give a decomposition of $\widehat{X}$ at certain stopping times which, following Schmidli [11], we call modified ladder epochs.

Let $\mathcal{P}(\mathcal{F})$ be the predictable $\sigma$-algebra on $\mathbb{R}_+ \times \Omega$ with respect to the filtration $\mathcal{F}$ introduced in Section 2. Let $\mathcal{B}_\partial$ be the Borel $\sigma$-algebra on $(0, \infty) \cup \{\partial\}$. If $\mathcal{H}: \mathbb{R}_+ \times \Omega \times ((0, \infty) \cup \{\partial\}) \to \mathbb{R}_+$ is a nonnegative process measurable with respect to $\mathcal{P}(\mathcal{F}) \otimes \mathcal{B}_\partial$, then the following compensation formula is valid (e.g., [8], page 439, or [2], page 9):

(4.2)
$$\mathbb{E}\bigg(\sum_{0 \leq t < \infty} \mathcal{H}(t, \omega, \Delta C(t, \omega))\bigg)$$
$$= \mathbb{E}\bigg(\int_0^\infty dt \int_{(0,\infty)} \nu(d\varepsilon) \mathcal{H}(t, \omega, \varepsilon)\bigg).$$

The first use of this formula will be to compute the expected number of times the new supremum of $\widehat{X}$ is attained by a jump of a subordinator $C$ over the previous supremum. Note that this is the case if and only if $\Delta C(t) > \widehat{S}(t-) - \widehat{X}(t-)$.

THEOREM 4.1. *The following formula is valid:*

$$(4.3) \qquad \mathbb{E}\left(\sum_{0 \leq t < \infty} \mathbb{1}_{\{\Delta C(t) > \widehat{S}(t-) - \widehat{X}(t-)\}}\right) = \frac{\mathbb{E}C(1)}{c - \mathbb{E}C(1)}.$$

PROOF. Take $\mathcal{H}(t, \omega, \varepsilon) := \mathbb{1}_{(\widehat{S}(t-,\omega) - \widehat{X}(t-,\omega), \infty)}(\varepsilon)$ in the compensation formula. The left-hand side in (4.2) is then precisely the left-hand side in (4.3). For the right-hand side in the compensation formula, compute

$$(4.4) \quad \begin{aligned} &\mathbb{E}\left(\int_0^\infty dt \int_{(0,\infty)} \nu(d\varepsilon) \mathbb{1}_{(\widehat{S}(t-) - \widehat{X}(t-), \infty)}(\varepsilon)\right) \\ &= \mathbb{E}\left(\int_0^\infty dt\, \nu(\widehat{S}(t-) - \widehat{X}(t-), \infty)\right) \\ &= \int_0^\infty \mathbb{E}[\nu(\widehat{S}(t) - \widehat{X}(t), \infty)]\, dt \\ &= \int_0^\infty dt\, \mathbb{E}[\nu(S(t), \infty)], \end{aligned}$$

where the third line follows by continuity in probability of $\widehat{X}$, and the fourth line by (4.1). Clearly, the last expression is equal to the monotone limit

$$(4.5) \quad \begin{aligned} &\lim_{q \to 0} \int_0^\infty e^{-qt}\, dt\, \mathbb{E}[\nu(S(t), \infty)] \\ &= \lim_{q \to 0} \frac{1}{q} \int_0^\infty qe^{-qt}\, dt\, \mathbb{E}[\nu(S(t), \infty)]. \end{aligned}$$

Let $\tau(q)$ be an exponential time with parameter $q$ independent of $C$ and $Z$, and let $F$ denote the distribution function of $S(\tau(q))$. Then $F$ is exponential with parameter $\Phi(q)$. It follows that

$$(4.6) \quad \begin{aligned} &\int_0^\infty qe^{-qt}\, dt\, \mathbb{E}[\nu(S(t), \infty)] \\ &= \mathbb{E}\left[\int_0^\infty qe^{-qt}\, dt\, \nu(S(t), \infty)\right] \\ &= \mathbb{E}[\nu(S(\tau(q)), \infty)] \\ &= \int_0^\infty \nu(x, \infty) F(dx) \\ &= \int_0^\infty (1 - e^{-\Phi(q)x}) \nu(dx), \end{aligned}$$

where the last equation follows by integration by parts. Further,

$$\lim_{q \to 0} \frac{1 - e^{-\Phi(q)x}}{q} = \lim_{q \to 0} \frac{\Phi(q)}{q} x = \frac{1}{\psi'(0+)} x = \frac{x}{d}.$$

By the monotone convergence theorem,

$$\lim_{q \to 0} \int_0^\infty e^{-qt} \, dt \, \mathbb{E}[\nu(S(t), \infty)]$$

(4.7)
$$= \lim_{q \to 0} \int_0^\infty \frac{1 - e^{-\Phi(q)x}}{q} \nu(dx)$$
$$= \frac{1}{d} \int_0^\infty x \nu(dx)$$
$$= \frac{\mathbb{E}C(1)}{c - \mathbb{E}C(1)}.$$

This proves formula (4.3). □

REMARK 4.2. We would like to emphasize a very important and somewhat subtle point which is a consequence of Theorem 4.1. Namely, the epochs when a new supremum of $\widehat{X}$ is reached by a jump of $C$ are discrete, and, in particular, neither time 0 nor any other time is an accumulation point of those epochs. More precisely, let us define

(4.8) $$\sigma_1 = \sigma := \inf\{t > 0 : \Delta C(t) > \widehat{S}(t-) - \widehat{X}(t-)\},$$

and inductively,

(4.9) $$\sigma_{n+1} := \inf\{t > \sigma_n : \Delta C(t) > \widehat{S}(t-) - \widehat{X}(t-)\}$$

on $\{\sigma_n < \infty\}$. Theorem 4.1 implies that $\sigma > 0$ a.s. and $\sigma_n < \sigma_{n+1}$ a.s. on $\{\sigma_n < \infty\}$. As a consequence, we can order the epochs when a new supremum is reached by a jump of a subordinator as $0 < \sigma_1 < \sigma_2 < \cdots$ a.s. The decomposition of $\widehat{X}$ and the ensuing derivations will depend on this result in a crucial way.

For $y > 0$, let $\widehat{\tau}_y := \inf\{t > 0 : \widehat{X}(t) > y\}$ be the entrance time of $\widehat{X}$ in $(y, \infty)$, and, similarly, $\tau_y := \inf\{t > 0 : X(t) > y\}$. Note that $\widehat{S}(t-) \leq y$ if and only if $t \leq \widehat{\tau}_y$. We need the expected occupation time formula for the reflected process $\widehat{S} - \widehat{X}$ before $\sigma \wedge \widehat{\tau}_y$.

PROPOSITION 4.3. *For $x > 0$ and $y > 0$, the following formula is valid:*

(4.10) $$\mathbb{E} \int_0^{\sigma \wedge \widehat{\tau}_y} \mathbb{1}_{(\widehat{S}(t) - \widehat{X}(t) \leq x)} \, dt = \mathbb{P}(\sigma = \infty, \widehat{\tau}_y = \infty) \frac{x}{d}.$$

PROOF. We first compute the expected occupation time of $\widehat{S} - \widehat{X}$ below $x$:

$$\mathbb{E} \int_0^\infty \mathbb{1}_{(\widehat{S}(t) - \widehat{X}(t) \leq x)} \, dt = \int_0^\infty \mathbb{P}(\widehat{S}(t) - \widehat{X}(t) \leq x) \, dt$$

(4.11)
$$= \int_0^\infty \mathbb{P}(S(t) \leq x) \, dt$$

$$= \mathbb{E}\tau_x.$$

Since $(\tau_x, x > 0)$ is a subordinator with the Laplace exponent $\Phi$, it follows that $\mathbb{E}\tau_x = (\mathbb{E}\tau_1)x = \Phi'(0+)x = x/d$.

Now we compute the expected occupation time of $\widehat{S} - \widehat{X}$ below $x$ after time $\sigma \vee \widehat{\tau}_y$:

$$\mathbb{E} \int_0^\infty \mathbb{1}_{(\widehat{S}(t) - \widehat{X}(t) \leq x)} \mathbb{1}_{(t > \sigma)} \mathbb{1}_{(\widehat{S}(t) > y)} \, dt$$

$$= \mathbb{E}\left[ \int_{\sigma \vee \widehat{\tau}_y}^\infty \mathbb{1}_{(\widehat{S}(t) - \widehat{X}(t) \leq x)} \, dt, \sigma \vee \widehat{\tau}_y < \infty \right]$$

(4.12)
$$= \mathbb{P}(\sigma \vee \widehat{\tau}_y < \infty) \mathbb{E}\left[ \int_{\sigma \vee \widehat{\tau}_y}^\infty \mathbb{1}_{(\widehat{S}(t) - \widehat{X}(t) \leq x)} \, dt | \sigma \vee \widehat{\tau}_y < \infty \right]$$

$$= \mathbb{P}(\sigma \vee \widehat{\tau}_y < \infty) \mathbb{E} \int_0^\infty \mathbb{1}_{(\widehat{S}(t) - \widehat{X}(t) \leq x)} \, dt$$

$$= \mathbb{P}(\sigma < \infty, \widehat{\tau}_y < \infty) \frac{x}{d}.$$

To justify the passage from the third to the fourth line, note that $\sigma \vee \widehat{\tau}_y$ is a stopping time at which $\widehat{X}$ reaches a new maximum, and hence by the strong Markov property, the reflected process $\widehat{S} - \widehat{X}$ starts afresh from 0. Similarly,

(4.13) $$\mathbb{E} \int_0^\infty \mathbb{1}_{(\widehat{S}(t) - \widehat{X}(t) \leq x)} \mathbb{1}_{(\widehat{S}(t) > y)} \, dt = \mathbb{P}(\widehat{\tau}_y < \infty) \frac{x}{d}.$$

Subtracting (4.12) from (4.13), it follows that

(4.14) $$\mathbb{E} \int_0^\infty \mathbb{1}_{(\widehat{S}(t) - \widehat{X}(t) \leq x)} \mathbb{1}_{(t \leq \sigma)} \mathbb{1}_{(\widehat{S}(t) > y)} \, dt = \mathbb{P}(\sigma = \infty, \widehat{\tau}_y < \infty) \frac{x}{d}.$$

One can prove similarly that

(4.15) $$\mathbb{E} \int_0^\infty \mathbb{1}_{(\widehat{S}(t) - \widehat{X}(t) \leq x)} \mathbb{1}_{(t \leq \sigma)} \, dt = \mathbb{P}(\sigma = \infty) \frac{x}{d}.$$

Finally, (4.10) follows by subtracting (4.14) from (4.15). □

Note that the proposition says that the expected occupation time measure for $\widehat{S} - \widehat{X}$ before $\sigma \wedge \widehat{\tau}_y$ is proportional to the Lebesgue measure on

$[0, \infty)$. Hence, formula (4.10) is by definition of the expected occupation time measure equivalent to

$$(4.16) \quad \mathbb{E} \int_0^{\sigma \wedge \widehat{\tau}_y} f(\widehat{S}(t) - \widehat{X}(t))\, dt = \frac{\mathbb{P}(\sigma = \infty, \widehat{\tau}_y = \infty)}{d} \int_0^\infty f(u)\, du,$$

where $f$ is a nonnegative Borel function on $[0, \infty)$.

Let $J := (\Delta C(\sigma) - (\widehat{S}(\sigma-) - \widehat{X}(\sigma-)))\mathbb{1}_{(\sigma < \infty)}$ be the overshoot at time $\sigma$. In the next proposition we compute the preliminary version of the joint distribution of the vector $(\widehat{S}(\sigma-), J, \widehat{S}(\sigma-) - \widehat{X}(\sigma-))$ on $\{\sigma < \infty\}$.

PROPOSITION 4.4. *For $x, y, z > 0$,*

$$(4.17) \quad \begin{aligned} &\mathbb{P}(\widehat{S}(\sigma-) \leq y, J > x, \widehat{S}(\sigma-) - \widehat{X}(\sigma-) > z, \sigma < \infty) \\ &= \frac{\mathbb{P}(\sigma = \infty, \widehat{\tau}_y = \infty)}{d} \int_{x+z}^\infty \nu(u, \infty)\, du. \end{aligned}$$

PROOF. We use the compensation formula with

$\mathcal{H}(t, \omega, \varepsilon)$

$:= \mathbb{1}_{(\widehat{S}(t-,\omega) \leq y)} \mathbb{1}_{(\widehat{S}(t-,\omega) - \widehat{X}(t-,\omega) > z)} \mathbb{1}_{(t \leq \sigma(\omega))} \mathbb{1}_{(x + \widehat{S}(t-,\omega) - \widehat{X}(t-,\omega), \infty)}(\varepsilon).$

Then

$$\mathbb{E} \sum_{0 \leq t < \infty} \mathcal{H}(t, \omega, \Delta C(t, \omega))$$
$$= \mathbb{P}(\widehat{S}(\sigma-) \leq y, \widehat{S}(\sigma-) - \widehat{X}(\sigma-) > z, J > x, \sigma < \infty).$$

On the other hand,

$$\begin{aligned} \mathbb{E}\bigg(\int_0^\infty dt &\int_{(0,\infty)} \nu(d\varepsilon) \mathcal{H}(t, \omega, \varepsilon)\bigg) \\ &= \mathbb{E}\bigg(\int_0^\sigma dt\, \mathbb{1}_{(\widehat{S}(t-) \leq y)} \mathbb{1}_{(\widehat{S}(t-) - \widehat{X}(t-) > z)} \\ &\quad \times \int_{(0,\infty)} \mathbb{1}_{(x + \widehat{S}(t-) - \widehat{X}(t-), \infty)}(\varepsilon) \nu(d\varepsilon)\bigg) \\ &= \mathbb{E}\bigg(\int_0^{\sigma \wedge \widehat{\tau}_y} dt\, \mathbb{1}_{(\widehat{S}(t) - \widehat{X}(t) > z)} \nu(x + \widehat{S}(t) - \widehat{X}(t), \infty)\bigg) \\ &= \frac{\mathbb{P}(\sigma = \infty, \widehat{\tau}_y = \infty)}{d} \int_0^\infty \mathbb{1}_{(z,\infty)}(u) \nu(x + u, \infty)\, du \\ &= \frac{\mathbb{P}(\sigma = \infty, \widehat{\tau}_y = \infty)}{d} \int_{x+z}^\infty \nu(u, \infty)\, du, \end{aligned}$$

where the fourth line follows from (4.16) with $f(u) = \mathbb{1}_{(z,\infty)}(u)\nu(x+u,\infty)$.
□

From formula (4.17) we can easily derive several useful corollaries.

COROLLARY 4.5. *The following formulae are valid:*

(4.18) $$\mathbb{P}(\sigma < \infty) = \rho,$$

(4.19) $$\mathbb{P}(J > x | \sigma < \infty) = \frac{1}{\mathbb{E}C(1)} \int_x^\infty \nu(u,\infty)\,du = 1 - H(x).$$

PROOF. Let $x \to 0$, $y \to \infty$ and $z \to 0$ in (4.17). It follows that
$$\mathbb{P}(\sigma < \infty) = \frac{\mathbb{P}(\sigma = \infty)}{d}\mathbb{E}C(1).$$
Solving for $\mathbb{P}(\sigma < \infty)$, we get (4.15). To obtain (4.19), let $y \to \infty$ and $z \to 0$ in (4.17). It follows that
$$\mathbb{P}(J > x, \sigma < \infty) = \frac{\mathbb{P}(\sigma = \infty)}{d}\int_x^\infty \nu(u,\infty)\,du.$$
By conditioning,
$$\mathbb{P}(J > x | \sigma < \infty) = \frac{1-\rho}{\rho d}\int_x^\infty \nu(u,\infty)\,du$$
$$= \frac{1}{\mathbb{E}C(1)}\int_x^\infty \nu(u,\infty)\,du.$$
□

In the next corollary, we interpret $\widehat{S}(\sigma-)$ as the absolute supremum $\widehat{S}(\infty)$ in case $\sigma = \infty$.

COROLLARY 4.6. *The event $\{\sigma < \infty\}$ and the random variable $\widehat{S}(\sigma-)$ are independent. As a consequence, the conditional distribution of $\widehat{S}(\sigma-)$ given $\sigma < \infty$ is equal to the unconditional distribution of $\widehat{S}(\sigma-)$.*

PROOF. Let $x \to 0$ and $z \to 0$ in (4.17). It follows that
$$\mathbb{P}(\widehat{S}(\sigma-) \leq y, \sigma < \infty)$$
(4.20)
$$= \mathbb{P}(\sigma = \infty, \widehat{\tau}_y = \infty)\frac{\mathbb{E}C(1)}{d}$$
$$= \mathbb{P}(\sigma = \infty, \widehat{S}(\infty) \leq y)\frac{\mathbb{E}C(1)}{d}.$$
Clearly,
$$\mathbb{P}(\widehat{S}(\sigma-) \leq y, \sigma = \infty) = \mathbb{P}(\widehat{S}(\infty) \leq y, \sigma = \infty).$$

Adding up,

$$\mathbb{P}(\widehat{S}(\sigma-) \leq y) = \left(\frac{\mathbb{E}C(1)}{d} + 1\right)\mathbb{P}(\widehat{S}(\infty) \leq y, \sigma = \infty)$$
$$= \frac{c}{d}\mathbb{P}(\widehat{S}(\infty) \leq y, \sigma = \infty).$$

Therefore,

$$\mathbb{P}(\widehat{S}(\sigma-) \leq y)\mathbb{P}(\sigma < \infty) = \frac{c}{d}\mathbb{P}(\widehat{S}(\infty) \leq y, \sigma = \infty)\frac{\mathbb{E}C(1)}{c}$$
$$= \frac{\mathbb{E}C(1)}{d}\mathbb{P}(\widehat{S}(\infty) \leq y, \sigma = \infty)$$
$$= \mathbb{P}(\widehat{S}(\sigma-) \leq y, \sigma < \infty)$$

by (4.20). $\square$

It follows that

(4.21)
$$\begin{aligned}\mathbb{P}(\sigma = \infty, \widehat{\tau}_y = \infty) &= \mathbb{P}(\sigma = \infty, \widehat{S}(\infty) \leq y)\\ &= \mathbb{P}(\sigma = \infty, \widehat{S}(\sigma-) \leq y)\\ &= \mathbb{P}(\sigma = \infty)\mathbb{P}(\widehat{S}(\sigma-) \leq y | \sigma = \infty)\\ &= \mathbb{P}(\sigma = \infty)\mathbb{P}(\widehat{S}(\sigma-) \leq y).\end{aligned}$$

Let $\widetilde{G}$ denote the distribution function of $\widehat{S}(\sigma-)$. Proposition 4.4 can be now improved as:

THEOREM 4.7. *The distribution of the vector $(\widehat{S}(\sigma-), J, \widehat{S}(\sigma-) - \widehat{X}(\sigma-))$ on the set $\{\sigma < \infty\}$ is given by*

(4.22)
$$\mathbb{P}(\widehat{S}(\sigma-) \leq y, J > x, \widehat{S}(\sigma-) - \widehat{X}(\sigma-) > z, \sigma < \infty)$$
$$= \mathbb{P}(\widehat{S}(\sigma-) \leq y)\left(\frac{1}{\mathbb{E}C(1)}\int_{x+z}^{\infty} \nu(u, \infty)\,du\right)\mathbb{P}(\sigma < \infty).$$

*Moreover, $\widehat{S}(\sigma-)$ and $J$ are conditionally independent given $\sigma < \infty$, and*

(4.23) $\quad\mathbb{P}(\widehat{S}(\sigma-) \leq y, J > x | \sigma < \infty) = \widetilde{G}(y)(1 - H(x)).$

REMARK 4.8. Formula (4.22) considerably extends the severity of ruin formula (see, e.g., [9], page 168).

It is now possible to write the absolute maximum of $\widehat{X}$ as a random sum of modified ladder heights. Recall that $\sigma_1 = \sigma$ and

$$\sigma_{n+1} = \inf\{t > \sigma_n : \Delta C(t) > \widehat{S}(t-) - \widehat{X}(t-)\}$$

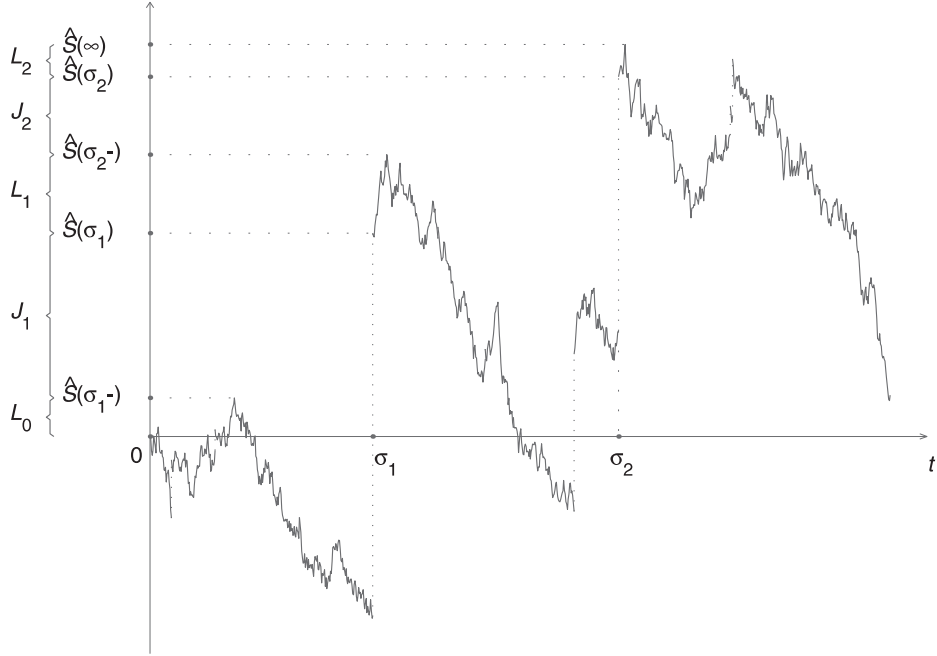

FIG. 1. *A trajectory of the dual process* $\widehat{X}(t) = -ct + C(t) + \varsigma W(t)$, *where* $C(t)$ *is a compound Poisson process,* $W(t)$ *is a standard Brownian motion and* $\varsigma > 0$.

on $\{\sigma_n < \infty\}$. Let $L_0 := \widehat{S}(\sigma_1-)$, $J_1 := \widehat{S}(\sigma_1) - \widehat{S}(\sigma_1-)$ and $L_1 := \widehat{S}(\sigma_2-) - \widehat{S}(\sigma_1)$ on $\{\sigma_1 < \infty\}$, and so on. We call $L_0, J_1, L_1, \ldots$ the modified ladder heights. Let also $N := \max\{n : \sigma_n < \infty\}$. By strong Markov property of $\widehat{X}$, $N$ has geometric distribution with parameter $\mathbb{P}(\sigma_1 = \infty) = 1 - \rho$. Clearly,

$$\widehat{S}(\infty) = L_0 + J_1 + L_1 + \cdots + J_N + L_N. \tag{4.24}$$

See Figure 1.

Note that $\mathbb{P}(L_0 \leq x, N = 0) = \mathbb{P}(\widehat{S}(\sigma-) \leq x, \sigma = \infty) = \widetilde{G}(x)(1-\rho)$. For every $n \in \mathbb{N}$, by the strong Markov property at $\sigma_n$, and by (4.23), we have

$$\mathbb{P}(L_0 + J_1 + L_1 + \cdots + J_n + L_n \leq x, N = n)$$
$$= (1-\rho)\rho^n (\widetilde{G}^{(n+1)*} * H^{n*})(x).$$

This leads to the Pollaczek–Hinchin formula for the distribution function of $\widehat{S}(\infty)$.

THEOREM 4.9. *For* $x \geq 0$,

$$\mathbb{P}(\widehat{S}(\infty) \leq x) = (1-\rho) \sum_{n=0}^{\infty} \rho^n (\widetilde{G}^{(n+1)*} * H^{n*})(x). \tag{4.25}$$

Equating this formula for $\theta(x) = \mathbb{P}(\widehat{S}(\infty) \leq x)$ and (3.4), one obtains

$$(1-\rho)\sum_{n=0}^{\infty}\rho^n(G^{(n+1)*} * H^{n*})(x) = (1-\rho)\sum_{n=0}^{\infty}\rho^n(\widetilde{G}^{(n+1)*} * H^{n*})(x).$$

By computing Laplace transforms of both sides, we get that

$$(4.26) \qquad \frac{(1-\rho)\mathcal{L}G(\beta)}{1-\rho\mathcal{L}G(\beta)\mathcal{L}H(\beta)} = \frac{(1-\rho)\mathcal{L}\widetilde{G}(\beta)}{1-\rho\mathcal{L}\widetilde{G}(\beta)\mathcal{L}H(\beta)}, \qquad \beta > 0,$$

from which it immediately follows that $\widetilde{G} = G$. Thus we have proved the following

COROLLARY 4.10. *The random variables*

$$\sup_{0\leq t<\infty}(-ct - Z(t)) \quad and \quad \sup_{0\leq t<\sigma}(-ct + C(t) - Z(t))$$

*have equal distributions.*

**5. Decomposition of the ladder height process.** In the previous section we looked at the process $\widehat{X}$ at the modified ladder epochs $\sigma_1 < \sigma_2 < \cdots$, and essentially decomposed the $\widehat{X}$ at these epochs. In this section we consider the ladder height process $\widehat{H}$ of $\widehat{X}$ obtained by time-changing the supremum process $\widehat{S}$ by the inverse local time at zero of the reflected process $\widehat{S} - \widehat{X}$. The excursion representation of the process $\widehat{S} - \widehat{X}$ will be combined with fluctuation identities and results from Section 4 to give a detailed description of the ladder height process.

Let us first briefly recall that $\widehat{S} - \widehat{X}$ is a strong Markov process, and hence it admits a local time process at 0, $\widehat{L} = (\widehat{L}(t), t \geq 0)$. The process $\widehat{L}$ is continuous and nondecreasing, and increases only when $\widehat{S} - \widehat{X}$ is at zero, or in other words, when $\widehat{X}$ reaches a new supremum. More precisely, the support of the measure $\widehat{L}(dt,\omega)$ is the zero set of $\widehat{S}(t,\omega) - \widehat{X}(t,\omega)$. If $\widehat{L}^{-1}(t) := \inf\{s > 0 : \widehat{L}(s) > t\}$ is the inverse of $\widehat{L}$ and one defines

$$(5.1) \qquad \widehat{H}(t) := \begin{cases} \widehat{S}(\widehat{L}^{-1}(t)), & \widehat{L}^{-1}(t) < \infty, \\ +\infty, & \text{otherwise,} \end{cases}$$

it is well known that the process $((\widehat{L}^{-1}(t), \widehat{H}(t)), t < \widehat{L}(\infty))$ is a two-dimensional subordinator killed at rate $q := 1/\mathbb{E}\widehat{L}(\infty)$ ([2], page 156). In particular, $\widehat{H} = (\widehat{H}(t), 0 \leq t < \widehat{L}(\infty))$ is a subordinator killed at rate $q = 1/\mathbb{E}\widehat{L}(\infty)$. Clearly,

$$(5.2) \qquad \widehat{H}(\widehat{L}(\infty)-) = \sup_{0\leq t<\widehat{L}(\infty)}\widehat{H}(t) = \sup_{0\leq t<\infty}\widehat{S}(t) = \widehat{S}(\infty),$$

and hence, the distribution function of $\widehat{S}(\infty)$ is equal to the distribution function of $\widehat{H}(\widehat{L}(\infty)-)$. Fluctuation identities give a formula for the Laplace exponent of $\widehat{H}$.

LEMMA 5.1. *The Laplace exponent $\widehat{\kappa}$ of $\widehat{H} = (\widehat{H}(t), 0 \leq t < \widehat{L}(\infty))$, with $\widehat{L}$ suitably normalized, is given by the following formula:*

$$(5.3) \quad \widehat{\kappa}(\beta) = \frac{\psi(\beta)}{\beta} = d + \mathbb{E}C(1) \int_{(0,\infty)} (1 - e^{-\beta x}) H(dx) + \frac{\psi_Z(\beta)}{\beta},$$

*where $H$ is the finite measure defined in (2.1).*

PROOF. The bivariate Laplace exponent $\widehat{\kappa}(\alpha, \beta)$ of $((\widehat{L}^{-1}(t), \widehat{H}(t)), t < \widehat{L}(\infty))$ is defined by

$$\exp\{-\widehat{\kappa}(\alpha, \beta)\} = \mathbb{E}[\exp -\{\alpha \widehat{L}^{-1}(1) + \beta \widehat{H}(1)\}], \qquad \alpha, \beta > 0.$$

The explicit formula for $\widehat{\kappa}$ comes from fluctuation theory:

$$(5.4) \qquad \widehat{\kappa}(\alpha, \beta) = k \frac{\alpha - \psi(\beta)}{\Phi(\alpha) - \beta},$$

where $k$ is a constant depending on the normalization of the local time. We take $k = 1$. By letting $\alpha = 0$ in (5.4), we obtain the Laplace exponent of $\widehat{H}$:

$$(5.5) \quad \begin{aligned} \widehat{\kappa}(\beta) &= \widehat{\kappa}(0, \beta) = \frac{\psi(\beta)}{\beta} \\ &= c - \frac{\Phi_C(\beta)}{\beta} + \frac{\psi_Z(\beta)}{\beta}, \qquad \beta > 0. \end{aligned}$$

Further, integrating by parts, we get

$$\begin{aligned} c - \frac{\Phi_C(\beta)}{\beta} &= c + \int_0^\infty (1 - e^{-\beta x}) \nu(x, \infty) \, dx - \int_0^\infty \nu(x, \infty) \, dx \\ &= c + \mathbb{E}C(1) \int_0^\infty (1 - e^{-\beta x}) H(dx) - \mathbb{E}C(1). \end{aligned}$$

Together with (5.5), this gives (5.3). □

REMARK 5.2. Note that the same argument shows that the ladder height process of $(-Z(t), t \geq 0)$ has the Laplace exponent equal to $\psi_Z(\beta)/\beta$.

Let $\widehat{L}$ be the local time of $\widehat{S} - \widehat{X}$ normalized by the choice $k = 1$ in formula (5.4). The excursion process $(e_s : s > 0)$ of the reflected process $\widehat{S} - \widehat{X}$ can be viewed as the superposition of three independent Poisson point processes: finite duration excursions that end with a jump of the subordinator, finite duration excursions that do not end with a jump of the subordinator, and excursions of infinite duration. Note that one can include the jump, if any, that concludes an excursion as part of that excursion and retain a Poisson point process. Also one needs to "carry along" the information about which jumps come from $C$ and which from $Z$, but that is a question of choosing a suitable filtration.

The excursion process and with it the ladder height process is killed at the time of the arrival of the first excursion of infinite duration. From the excursion picture we know that $\widehat{L}(\infty)$ is an exponential random variable with parameter equal to the killing rate of the ladder height process. By convention, $\widehat{H}(t) = \infty$ after killing. It is easily shown that the killing rate can be obtained from the Laplace exponent as $\psi(0+)$. Since $\mathbb{E}Z(t) = 0$, we know that $\psi_Z(\beta)/\beta \to 0$ as $\beta \to 0$. The measure $H$ is finite, so by dominated convergence,

$$\int_{(0,\infty)} (1 - e^{-\beta x}) H(dx) \to 0$$

as $\beta \to 0$. Using the explicit formula (5.3) one finds that the killing rate equals $d = c - \mathbb{E}C(1)$, or in other words, $\widehat{L}(\infty) \sim \exp(c - \mathbb{E}C(1))$.

Recall that $\sigma$ is the first modified ladder epoch. On the local time scale $\widehat{L}(\sigma)$ corresponds to the time of the first arrival of an excursion that ends with a jump of the subordinator $C$ unless the excursion process is killed first. The probability $\mathbb{P}(\sigma < \infty)$ is therefore equal to the probability that the first excursion that ends with a jump of the subordinator arrives before the first excursion of infinite duration. We are thus computing the probability that one of the two independent Poisson processes will "claim" the first arrival. It is well known that we can compute

$$\rho = \mathbb{P}(\sigma < \infty) = \mathbb{P}(\widehat{L}(\sigma) < \widehat{L}(\infty)),$$

where $\widehat{L}(\sigma)$ and $\widehat{L}(\infty)$ are two independent exponential random variables. Since we know the rate of $\widehat{L}(\infty)$, a simple calculation shows that the rate of arrival of excursions that end with a jump of $C$ equals $\mathbb{E}C(1)$.

Let us turn to the ladder height process $\widehat{H}$. The jumps of $\widehat{H}$ are a Poisson point process. These jumps, however, can be seen as a mapping of the "bigger" process $e$ of excursions. By the mapping theorem (see [6], page 17) the jumps of $\widehat{H}$ are a superposition of two independent Poisson processes killed at an independent exponential time of rate $d = c - \mathbb{E}C(1)$. The mapping theorem is applied in such a way that the image of $e_s$ is $\partial$ if the excursion ends without a jump. The jumps coming from excursions that end with a jump of $C$ contribute a pure jump process to $\widehat{H}$. The other excursions contribute jumps and possibly a deterministic drift. But the jumps coming from $Z$ are an independent Poisson process. Hence the process $\widehat{H}$ is a sum of two independent subordinators $\eta$ and $\zeta$ killed at an independent exponential time $\tau := \widehat{L}(\infty)$.

The subordinator $\eta$ corresponds to increases of $\widehat{S}$ due to jumps of $C$. As the set of times when $\widehat{S}$ increases by a jump of $C$ is discrete, the process $\eta$ is compound Poisson with arrival rate equal to $\mathbb{E}C(1)$. Jumps of $\eta$ have the same distribution as the overshoot $J$ which is given by the formula (4.19): $\mathbb{P}(J \in dx | \sigma < \infty) = \nu(x, \infty)\, dx / \mathbb{E}C(1)$. This means that $\eta$ contributes exactly

$$\mathbb{E}C(1) \int_{(0,\infty)} (1 - e^{-\beta x}) H(dx)$$

to the Laplace exponent of the ladder height process $\widehat{H}$ given in (5.3).

The subordinator $\zeta$ arising from increases of $\widehat{H}$ not due to jumps of $C$ is independent of $\eta$. This leaves us with the conclusion, given that the killing rate is $c - \mathbb{E}C(1)$, that the Laplace exponent of $\zeta$ is $\psi_Z(\beta)/\beta$. This way the groundwork in Section 4 has been translated into a decomposition of the ladder height process:

THEOREM 5.3. *Let $\widehat{X}(t) = -ct + C(t) - Z(t)$ and let $\widehat{H}$ be the ladder height process of $\widehat{X}$. The following assertions are valid:*

(i) *$\widehat{H}$ is killed at rate $d = c - \mathbb{E}C(1)$.*

(ii) *The ladder height process is the sum of two independent subordinators $(\widehat{H}(t) = \eta(t) + \zeta(t) : 0 \leq t < \widehat{L}(\infty))$. The subordinator $\eta$ corresponds to jumps of $\widehat{S}$ due to the claim process $C$. It is compound Poisson with Lévy measure $\mathbb{E}C(1)H$. The subordinator $\zeta$ corresponds to increases of $\widehat{S}$ not arising from jumps of $C$. Its Laplace exponent is $\psi_Z(\beta)/\beta$.*

The decomposition of $\widehat{H}$ gives some insights into the structure of the process $\widehat{H}$. From the form of the Laplace exponent of $\zeta$, we find that the distribution of $\zeta$ does not depend in any way on the distribution of $C$, which is a remarkable conclusion given that $C$ is a subordinator with a dense set of times of jumps. Many other conclusions from Section 4 can be recast in terms of the ladder height process.

As an example, consider Corollary 4.6. Recall that $\tau = \widehat{L}(\infty)$. Define $\gamma := \widehat{L}(\sigma)$ if $\sigma < \infty$ and $\infty$ else. The event $\{\sigma < \infty\}$ is equal to $\{\gamma < \tau\}$. It is easily seen that

$$\widehat{S}(\sigma-) = \zeta(\gamma \wedge \tau).$$

Note that $\widehat{S}(\sigma-) = \widehat{S}(\infty)$ on $\{\sigma = \infty\}$. Combining Corollary 4.6 and the conclusion following (4.26), we know that the distribution function of $\widehat{S}(\sigma-)$ is $G$. The decomposition of $\widehat{H}$ gives further information about $G$. Since the random variables $\gamma$ and $\tau$ are independent of $\zeta$, and $\gamma \wedge \tau \sim \exp(c)$, the distribution of $\widehat{S}(\sigma-)$ is that of a subordinator taken at an independent exponential time and hence infinitely divisible with Lévy measure $\Lambda$ given

by

$$(5.6) \qquad \Lambda(dx) = \int_{(0,\infty)} \frac{e^{-ct}}{t} \mathbb{P}(\zeta_t \in dx)\, dt.$$

See [2], page 162. However, a direct computation using the Laplace exponent of $\zeta$ yields

$$(5.7) \qquad \mathbb{E}(\exp(-\beta \zeta_{\gamma \wedge \tau})) = \int_{(0,\infty)} \mathbb{E}(\exp(-\beta \zeta_t)) c e^{-ct}\, dt$$

$$= c \int_{(0,\infty)} \exp(-t\psi_Z(\beta)/\beta) e^{-ct}\, dt$$

$$(5.8) \qquad = \frac{c\beta}{c\beta + \psi_Z(\beta)}.$$

Comparing this formula to (3.2) gives an independent proof of Corollary 4.10.

M. Huzak  
H. Šikić  
Z. Vondraček  
Department of Mathematics  
University of Zagreb  
Bijenička c. 30  
10000 Zagreb  
Croatia  
e-mail: huzak@math.hr  
e-mail: hsikic@math.hr  
e-mail: vondra@math.hr

M. Perman  
Institute for Mathematics,  
  Physics and Mechanics  
University of Ljubljana  
Jadranska 19  
1000 Ljubljana  
Slovenia  
e-mail: mihael@valjhun.fmf.uni-lj.si